\documentclass[12pt]{article}

\providecommand{\keywords}[1]{\textbf{\textit{Keywords---}} #1}

\usepackage{amssymb}
\usepackage{amsfonts}
\usepackage{amsmath}
\usepackage{amssymb}
\usepackage{setspace}
\usepackage{graphicx}

\usepackage{amsmath, amssymb, amsthm, latexsym, mathrsfs, color}

\newtheorem{definition1}{Definition}[section]
\newtheorem{lemma1}[definition1]{Lemma}
\newtheorem{theorem1}[definition1]{Theorem}
\newtheorem{proposition1}[definition1]{Proposition}

\newtheorem{remark1}[definition1]{Remark}
\begin{document}

\title{Countable Partially Exchangeable Mixtures}

\author{Cecilia Prosdocimi \thanks{Dipartimento di Economia e Finanza, Universit\`{a} LUISS, viale Romania 32, 00197 Roma, Italy, 
              prosdocimi.cecilia@gmail.com}
                      \and
        Lorenzo Finesso \thanks{ Istituto di Elettronica e di Ingegneria dell'Informazione e delle Telecomunicazioni, Consiglio Nazionale delle Ricerche, via Gradenigo 6/a, Padova, Italy,
         lorenzo.finesso@ieiit.cnr.it}
} 
\date{\today}

\maketitle
%


\begin{abstract}
Partially exchangeable sequences representable as mixtures of Markov
chains are completely specified by de Finetti's mixing measure. The
paper characterizes, in terms of a subclass of hidden Markov models,
the partially exchangeable sequences with mixing measure
concentrated on a countable set,
for sequences of random variables both on a discrete and a
Polish space.
\end{abstract}


\section{Introduction}

In the Hewitt-Savage generalization of de Finetti's theorem,
the distributions of exchangeable sequences of Polish valued random variables are shown to be in one to one
correspondence with distributions of mixtures of i.i.d. sequences, and ultimately
with the mixing measure defining the mixture. The mixing measure
thus acts as a model of the exchangeable sequence and its
properties shed light on the random mechanism generating the
sequence. In this regard~\cite{DIACONIS-FREEDMAN-regularity1},
using the connection with the Markov moment problem, characterizes
the subclass of exchangeable sequences of discrete valued random variables whose mixing measures are
absolutely continuous, and have densities in $L^p$.
It is also of
interest to characterize the subclass for which the mixing measure
is discrete (\textit{i.e.} concentrated on a countable set). A contribution in this direction has been given
by~\cite{DHARMADHIKARI}, where it is proved that an
exchangeable sequence of discrete valued random variables has de Finetti mixing measure concentrated on
a countable set if and only if it is a hidden Markov model (HMM),
see Definition \ref{DEF_HMM} for the precise notion.

\medskip

The more general class of partially exchangeable sequences is in one to one
correspondence with mixtures of Markov chains.
As noted in~\cite{DIACONIS-FREEDMAN-regularity1}, the results on the regularity of the mixing measure
carry to partially exchangeable sequences. On the other hand, to the best of
our knowledge, no results have been reported concerning partially exchangeable sequences
with discrete mixing measures.

The goal of the present paper is to
characterize, in the spirit of~\cite{DHARMADHIKARI}, countable mixtures of Markov chains.
Our results hold for sequences of both discrete and Polish space valued random variables.
This has required
the development of a few special results for HMMs, previously not available in the literature. Of
independent interest are Propositions~\ref{LEMMA_SUCCESSIVES MATRIX} and~\ref{PROPOSITION_SUCCESSIVES MATRIX_CONTINUOUS_STATE_SPACE}
 on the rows of the array of the successors of an HMM,
and most of the Appendix, on properties of sequences of stopping times with respect to the filtration generated by
an HMM and its underlying Markov chain.
For sequences of Polish valued random variables it has also been necessary to first extend~\cite{DHARMADHIKARI}
to show the equivalence between exchangeable HMMs and countable mixtures of i.i.d. sequences.

\medskip
In Section~\ref{tools} we review the basic definitions in the setup most
convenient for our purpose. The reader should be aware of the fact
that slightly different notions of partial exchangeability coexist
in the literature for sequences of discrete valued random variables.
We recall the original definition, introduced in~\cite{DEFINETTI_1938} and elaborated
in~\cite{FORTINI_LADELLI_PETRIS_REGAZZINI}. The latter paper
clarifies the relationship with the alternative definition given
in~\cite{DIACONIS_FREEDMAN_1980}. Section \ref{Countable state space}
deals with sequences of discrete valued random variables, Section \ref{Polish state spaces} with sequences
of Polish valued random variables.
In Section~\ref{Hidden Markov
Models} we constructively prove Proposition~\ref{LEMMA_SUCCESSIVES MATRIX}, which
is instrumental in the balance of the paper. In Section
\ref{Representations of countable Markov mixtures} we
characterize countable mixtures of Markov chains taking values on a discrete space.
In Section \ref{Polish state spaces}
we extend the result of~\cite{DHARMADHIKARI} to
exchangeable sequences of Polish valued random variables.
This allows  us to prove Theorem \ref{THEOREM_MAIN_CONTINUOUS_STATE_SPACE},
the main result of the paper.
Section~\ref{sec:conclus} contains final remarks and
hints at possible extensions.
In the Appendix are collected results on sequences of stopping times for HMMs,
unavailable elsewhere in the literature,
and a technical result for a special class of HMMs representable as mixtures of i.i.d. sequences.

\section{Preliminaries} \label{tools}
\setcounter{equation}{0}

$Y=(Y_n)_{n\ge 0}$ denotes a sequence of random variables on a
probability space $(\Omega,\mathcal F, \mathbb P )$, taking values
in a Polish space $S$ endowed with the Borel $\sigma$-field $\mathcal S$. The
generic element of $S$ is denoted $y$, and $y_0^N= y_0 y_1
\dots y_N$ is an element (a string) of the $N+1$-th fold Cartesian
product $S^{N+1}$, likewise $(Y_0^N = y_0^N)$ is the event
$(Y_0=y_0,Y_1=y_1,\dots,Y_N=y_N)$.

\medskip\noindent
\textit{Exchangeable and partially exchangeable sequences.} The notions of exchangeability and partial exchangeability both originate in the work of de Finetti. The former notion is well established in the literature, the latter has been presented in various disguises. For ease of reference we give below the definitions used in the paper.

\begin{definition1}\label{EXCHANGEABLE_SEQUENCES_GENERAL_STATE_SPACE}[\cite{KALLENBERG_BOOK_2005}, page 24]
A sequence of random variables $(Y_n)$ on $\big(S,\mathcal S\big)$
is \emph{exchangeable} if, for any $N \in \mathbb N$, and distinct natural numbers $m_1, \dots, m_N$,
\begin{equation*}
(Y_{m_1}, \dots, Y_{m_N})\stackrel{d} {=} (Y_1, \dots, Y_N),
\end{equation*}
where $\stackrel{d}{=}$ denotes equality in distribution.
\end{definition1}

We adopt the definition of partial exchangeability given in
\cite{FORTINI_LADELLI_PETRIS_REGAZZINI}, close to de Finetti's original \cite{DEFINETTI_1938}. The relation with the definition of partial exchangeability given in \cite{DIACONIS_FREEDMAN_1980} is clarified in
\cite{FORTINI_LADELLI_PETRIS_REGAZZINI}.
The definition requires the introduction of the \emph{successors array} $V$  of a given random sequence $(Y_n)$ on $S$,
and the extension of $S$ to $S^* = S \cup\{\partial\}$,
where $\partial \notin S$ is a fictitious state.

If $S$ is discrete (as in Section \ref{Countable state space} of the paper)
$V=(V_{y, n})_{y \in S, n \geq 1}$ is defined setting $V_{y, n}$
equal to the value of $Y$ immediately following the $n$-th visit to $y$ of $Y$.
If $Y$ visits $y$ only $m <\infty$ time,
to avoid rows of $V$ of finite length, one assigns
$V_{y, n}= \partial$ for all $n>m$ . 

If $S$ is uncountable (as in Section \ref{Polish state spaces} of the paper)
let $\mathcal E =\{E_j\}_{j \geq0}$ be a fixed countable partition of $S^*$,
with $E_0= \{  \partial\} $.
The matrix $V=(V_{j, n})_{j,n\geq1 }$ is defined setting $V_{j, n}$
equal to the value of $Y$ immediately following the $n$-th visit to $E_j$.
As in the discrete case, if $Y$ visits $E_j$ only $m <\infty$ times,
 set $V_{j, n}= \partial$ for all $n>m$.
In the uncountable case $V$ depends on the partition $\mathcal E$.

\begin{definition1} \label{DEF_PARTIAL_EXCHANGEABILITY} [\cite{FORTINI_LADELLI_PETRIS_REGAZZINI}]
A sequence of random variables $(Y_n)$ on $\big(S,\mathcal S\big)$ is \emph{partially exchangeable} if its successors array $V$
is distributionally invariant under finite, not necessarily identical,  permutations within each of its rows.
\end{definition1}

Note that, if $(Y_n)$ is partially exchangeable, the rows of its
successors array $V$ are exchangeable.

\medskip\noindent
\textit{Mixtures of i.i.d. sequences and of Markov chains.}
The set of all probability measures on $(S,\mathcal S)$ is denoted $\mathcal M_{\mathcal S}$
and it is equipped with the $\sigma$-field generated by the maps $p
\mapsto p(A)$, varying $p$ in $\mathcal M_{\mathcal S}$ and $A$ in $\mathcal S$.
\begin{definition1}\label{DEF_MIXTURES_IID_CONTINUOUS STATE SPACE}
$(Y_n)$ is a \textit{mixture of i.i.d. sequences} if there exists a random probability measure
 $\widetilde p$
on $\mathcal M_{\mathcal S}$ such that, for any $N\in \mathbb N$ and any $A_0, \dots, A_N  \in \mathcal S$,
\begin{eqnarray} \label{EQ_IID_MIXTURE}
\mathbb P \big (Y_0 \in A_0, \dots, Y_N \in A_N \mid \widetilde p) =
\widetilde p (A_0 )\dots \widetilde p (A_N), \qquad \mathbb P-a.s.
\end{eqnarray}
\end{definition1}
\begin{definition1}\label{DEF_COUNTABLE_MIXTURES_IID_CONTINUOUS STATE SPACE}
A mixture of i.i.d. sequences is \emph{countable (finite)} if the
random probability measure $\widetilde p$ is concentrated on a countable (finite) subset of $\mathcal M_{\mathcal S}$.
\end{definition1}

If $(Y_n)$ is a countable (finite) mixture of i.i.d. sequences  let $\big (p_h(\cdot)\big)_{h \in H}$, where $H$ is countable (finite),
be the set of measures on which $\widetilde p$ is concentrated then, integrating Equation (\ref{EQ_IID_MIXTURE}) over $\Omega$, one has
\begin{eqnarray} \label{IID_COUNTABLE_MIXTURE_CONTINUOUS_STATE_SPACE}
\mathbb P \big (Y_0 \in A_0, \dots, Y_N \in A_N  ) = \sum_{h \in H}  \mu_h p_h(A_0) \dots  p_h(A_N)  ,
\end{eqnarray}
where $\mu_h := \mathbb P (\widetilde p = p_h)>0$ and $\sum_{h \in H}\mu_h=1$.

\medskip
A time homogenous Markov chain with values in $S^*$ is characterized by a transition kernel
$k: S^{*} \times \mathcal S^{*}  \rightarrow [0,1]$.
For the purposes of the paper it is enough to consider Markov chains whose transition
kernels are constant, with respect to the first variable, on the elements of a fixed partition $\mathcal E=(E_j)_{j \geq 0}$ of $S$.
Kernels $k$ in this subclass can be represented in terms of the simpler kernels in the class $T^*$ made up of all kernels
$t: \mathbb N_0 \times \mathcal S^* \rightarrow [0,1]$,
where $\mathbb N_0 = \mathbb N  \cup \{0\}$. For a given $t\in T^*$ one defines $k=k_t$ as follows
\begin{equation}\label{DEF_k}
k_t(y,A) := \sum_{j \geq 0 } \mathbb I_{E_j}(y)\, t(j,A),
\end{equation}
where $\mathbb I_{E_j}(\cdot)$ is the indicator function of $E_j$.
The reader is referred to \cite{FORTINI_LADELLI_PETRIS_REGAZZINI} for a more detailed discussion.
We are ready to give the following

\begin{definition1}\label{DEF_MCs_MIXTURE_RANDOM_ELEMENT}
Let $(Y_n)$ be a sequence of random variables with $\mathbb P(Y_0=y_0)=1$, for some $y_0 \in E_1$, then $(Y_n)$ is a
\textit{mixture of homogeneous Markov chains} if there exists a
random kernel $\widetilde t$
on $T^*$ such that, for any $N\in \mathbb N$, and $A_1, \dots, A_N \in \mathcal S^*$,
\begin{equation}\label{EQ_MIX_MCs}
\mathbb P ( Y_1\in A_1, \dots, Y_N\in A_N \mid \widetilde t) =
\int_{A_1}\cdot\cdot \int_{A_N} k_{\widetilde t}(y_0, dy_1)\cdots k_{\widetilde t}(y_{N-1}, dy_N) \quad
\mbox{$\mathbb P$-a.s.}
\end{equation}
\end{definition1}
By the definition of $k_{ t}$ in (\ref{DEF_k}), Equation (\ref{EQ_MIX_MCs})
gives
\begin{equation}\label{EQ_MIX_MCs_t}
\mathbb P ( Y_1\in A_1, \dots, Y_N\in A_N \mid \widetilde t) =
\sum_{j_1^N \in (\mathbb N_0)^N}\int_{A_1\cap E_{j_1}}\dots \int_{A_N \cap E_{j_N}} \widetilde t(1, dy_1)\cdots\widetilde t(j_{N-1}, dy_N) \quad
\mbox{$\mathbb P$-a.s.}
\end{equation}

\begin{definition1}
A mixture of Markov chains is \emph{countable (finite)} if the random kernel
$\widetilde t$ is
 countably (finitely) valued.
 \end{definition1}

For a countable (finite) mixture of Markov chains,
Equation (\ref{EQ_MIX_MCs}) reads,
after integration over $\Omega$,
\begin{eqnarray}\label{EQ_COUNTABLE_MARKOV_MIXTURES}
\begin{split}
\mathbb P ( Y_1 \in A_1, \dots, &Y_N \in A_N  )\\=
&\sum_{h \in H}  \mu_h
\int_{A_1}\cdots\int_{A_N}
  k_{ t_h}(y_0, dy_1)\cdots k_{ t_h}(y_{N-1}, dy_N)
\end{split}
\end{eqnarray}
where $\left(t_h\right)_{h \in H}$ are the kernels on which $\widetilde t$ is concentrated,
$\mu_h:= \mathbb P (\widetilde t =t_h)$,
and $H$ is a countable (finite) set.
Finally we can write the finite distributions of a
countable mixtures of Markov chains as
\begin{eqnarray}\label{EQ_COUNTABLE_MARKOV_MIXTURES_EXPLICIT}
\begin{split}
\mathbb P ( Y_1 \in A_1, &\dots,  Y_N \in A_N  )\\=
&\sum_{h \in H}  \mu_h
\sum_{j_1^N \in (\mathbb N_0)^N}\int_{A_1\cap E_{j_1}}\dots \int_{A_N \cap E_{j_N}}  t_h(1, dy_1)\cdots t_h(j_{N-1}, dy_N) .
\end{split}
\end{eqnarray}

\begin{remark1}\label{REMARK_COUNTABLE_STATE_SPACE_MIXTURE_MCS}
To define mixtures of Markov chains when $S^*$ is discrete, less technicalities are needed since transition kernels reduce to
transition matrices. When $S^*$ is discrete Equation (\ref{EQ_MIX_MCs}) reduces to
\begin{equation}
\mathbb P ( Y_1^N = y_1^N \mid \widetilde P) = \widetilde P_{y_0y_1}
\widetilde P_{y_1y_2}   \dots \widetilde P_{y_{N-1}y_N} \qquad
\mbox{$\mathbb P$-a.s.},
\end{equation}
where $\widetilde P= (\widetilde P_{i,j})_{i,j \in S^*}$ is a random matrix
varying on the set $\mathcal P^*$ of transition matrices, and $P(Y_0=y_0)=1$ for some $y_0\in S$.
The analog of Equation~(\ref{EQ_COUNTABLE_MARKOV_MIXTURES}) is
\begin{equation}\label{EQ_COUNTABLE_MARKOV_MIXTURES_COUNTABLE_STATE_SPACE}
\mathbb P ( Y_1^N = y_1^N ) =
\sum_{h \in H}  \mu_h
 P^h_{y_0y_1} \dots  P^h_{y_{N-1}y_N},
\end{equation}
where $(P^h)_{h \in H}$ are the possible values taken by $\widetilde P$.
\end{remark1}

\medskip

\noindent \textit{Representation theorems.}
The classic de Finetti's representation theorem characterizes exchangeable sequences
as mixtures of i.i.d. sequences.
The statement, the proof, and an extensive discussion of the ramifications of the theorem
can be found in~\cite{ALDOUS} and~\cite{KALLENBERG_BOOK_2005}.
The representation theorem characterizing partially exchangeable sequences as mixtures of Markov chains
was first proved in \cite{DIACONIS_FREEDMAN_1980}, for $S$ discrete and using a slightly different notion of partial exchangeability, and later extended to $S$ a general Polish space in \cite{FORTINI_LADELLI_PETRIS_REGAZZINI}.

\medskip

\noindent
\textit{Hidden Markov models.}
Let $S$ be a Polish space endowed with the Borel $\sigma$-field $\mathcal S$.
\begin{definition1}\label{DEF_HMM}
A random sequence $(Y_n)$ is an hidden Markov model (HMM)
if there exists a pair $(X_n,\widetilde Y_n)$ taking values on $\mathcal X \times S$
for a discrete space $\mathcal X$,
such that
\begin{enumerate}
\item \,\, $(X_n) $ is a time homogeneous Markov chain on $\mathcal X$,
\item \,\,(\textit{conditional independence property}) for any $N \in \mathbb N$, and
for any $S_0, \dots, S_N \in \mathcal S$ it holds
\begin{equation*}
\mathbb P \big(\widetilde Y_0 \in S_0, \dots ,\widetilde Y_N \in S_N \mid X_0^N = x_0^N \big)
= \prod_{n=0}^N  \mathbb P \big(\widetilde Y_n \in S_n \mid X_n=x_n\big),
\end{equation*}
\item \,\, $(Y_n)$ and $(\widetilde Y_n)$ have the same distributions, i.e.
 for any $N \in \mathbb N$, and
for any $S_0, \dots, S_N \in \mathcal S$ it holds
$$
P\big(Y_0 \in S_0, \dots Y_N \in S_N\big) = P\big(\widetilde Y_0 \in S_0, \dots \widetilde Y_N \in S_N\big).
$$
\end{enumerate}
\end{definition1}
It is often possible to verify the second property directly
for the sequence $(Y_n)$.


 A HMM is characterized by the initial distribution $\pi$ on $\mathcal X$,
by the transition matrix $P=(P_{ij})_{i,j \in \mathcal X}$ of the Markov chain $(X_n)$,
and by the \textit{read-out} distributions $f_x (\bar S)$, with $\bar S \in \mathcal S$,
where

\begin{equation}\label{DEF_READ_OUT_DISTRIBUTIONS}
f_x (\bar S) := \mathbb P \big(\widetilde Y_n \in \bar S \mid X_n=x\big) .
\end{equation}
\noindent
We refer to the sequence $(X_n)$ as the "underlying Markov
chain" of the HMM.

\medskip

For discrete $S$,
there are many equivalent definitions of HMMs (see \cite{VIDYASAGAR_2011_REALIZATION}),
not all making sense for $S$ Polish.

\begin{remark1}\label{RMK_HMM_MIXT_IID}
A countable mixture of i.i.d. sequences
as in Equation (\ref{IID_COUNTABLE_MIXTURE_CONTINUOUS_STATE_SPACE})
is an HMM: take as $(X_n)$
the Markov chain with values in $H$, with identity transition matrix,
initial distribution $(\mu_1, \dots, \mu_h,\dots)_{h \in H}$,
and read-out distributions $ \mathbb P \big(Y_n \in \bar S \mid X_n =h \big)= p_h(\bar S)$.
\end{remark1}


\section{Countable Markov mixtures with discrete state space}\label{Countable state space}
\setcounter{equation}{0}
In this section $S$ is a discrete set.

\subsection{The successors array of hidden Markov models}\label{Hidden Markov Models}

The following result will be instrumental later, and it is also of independent interest.
It is based on some useful properties of HMMs, that can be found in the Appendix.
HMMs and the succesors array are defined in Section \ref{tools}.
\begin{proposition1} \label{LEMMA_SUCCESSIVES MATRIX}
Let $(Y_n)$ be a HMM on a discrete space $S$ with recurrent
\footnote{A time homogeneous Markov chain $(X_n)$ is recurrent if $\mathbb P (X_n= x \, \,\mbox{i.o. $n$} \mid X_1=x)=1$,
for all $x \in \mathcal X$
such that $\mathbb P (X_1=x)>0$.
Such Markov chains have no transient
states but possibly more than one recurrence class.}
underlying Markov chain, then each row of the successors array $(V_{y,n})$
is a HMM with recurrent
underlying Markov chain.
\end{proposition1}

\begin{proof} 
Denote with $\mathcal X$ the discrete state space of the Markov chain $(X_n)$, underlying the process $(Y_n)$ and let,
for any $x \in \mathcal X$ and $ y \in S$,
\begin{equation*}
f_x(y) := \mathbb P \big(  Y_n =  y \mid X_n = x \big)
\end{equation*}
be the read-out distributions. Fix $y \in S$.
To prove the theorem we construct a recurrent Markov chain $(W^y_n)_{n \geq 1}$
such that the pair $\left(  W^y_n,V_{y,n}\right)_{n \geq 1}$ satisfies the conditions
of Definition \ref{DEF_HMM}
of HMM
(note that for convenience we let time start at $n=1$).
The proof is divided in three main steps.

\smallskip

\textbf{Step 1} \textit{Construction of the Markov chain $(W^y_{ n})$.}
To construct the Markov chain $(W^y_n)$, define inductively the random times of
the $n$-th visit of $(Y_n)$ to state $y$:
\begin{equation}
\nonumber \tau_1^y := \inf \{ t \ge 0\mid Y_t=y  \} ,
\end{equation}
\begin{equation}
\nonumber \tau_n^y := \inf\{ t > \tau_{n-1}^y \mid Y_t=y \},
\end{equation}
with the usual convention $\inf \emptyset = + \infty$. The random
times $(\tau_n^y)$ are stopping times with respect to the
filtration spanned by $(Y_n)$,
and so are the times $(\tau_n^y +1)$.
The random times $(\tau_{n}^y+1)$ are actually hitting-times of $ \mathcal X \times y $,
according to the Definition \ref{DEF_HITTING_TIMES} in the Appendix.
Define the sequence
\begin{eqnarray}
W^y_{ n} := \left\{
\begin{array}{rl}
\varepsilon & \text{for }   \tau_n^y = +\infty, \\
X_{\tau^y_n+1} & \text{for }  \tau_n^y < +\infty,
\end{array} \right.
\end{eqnarray}
where $\varepsilon \notin \mathcal X$ is a fictitious state.
The sequence $W^y_{n}$ is either identically equal to
$\varepsilon$, or it never hits it since the times $\tau_n^y$ are either all finite or all infinite
\footnote{If $Y_{n}=y$, for some finite $n$, then $X_{n}= x$, for some $x$ such that $f_x(y)>0$. Since $X$ is recurrent it hits $x$ infinitely many times, thus $Y$ hits $y$ infinitely many times.}.

Let us check that $(W^y_{ n})$ is a Markov chain.
If the case $W^y_{n}\equiv \varepsilon$ obtains, $(W^y_{n})$ is a (recurrent) Markov chain.
Otherwise a direct computation gives, for any $ N \in \mathbb N$,
and any $x_1,\dots, x_N \in \mathcal X$,
\begin{multline*}
\mathbb P\big(  W^y_{ N}=x_N \mid W^y_{N-1 }=x_{N-1},\, \dots, W^y_{1}=x_1 \big) \\
=
\mathbb P\big(  X_{\tau^y_N+1}=x_N \mid X_{\tau^y_{N-1}+1}=x_{N-1},\,  \dots ,X_{\tau^y_{1}+1}=x_1\big) \\
=\mathbb P\big(  X_{\tau^y_N+1}=x_N \mid X_{\tau^y_{N-1}+1}=x_{N-1} \big) = \mathbb P
\big( W^y_{ N}=x_N \mid W^y_{ N-1 }=x_{N-1} \big),
\end{multline*}
where Remark~\ref{REMARK_STRONG_STRONG_MARKOV} in the Appendix applies.
Thus $(W^y_{ n})$ is Markov.

\smallskip

\textbf{Step 2} \textit{Check of the recurrence of $(W^y_{ n})$.} Since
\begin{equation*}
\mathbb P \left( W^y_{ n}=x \,\,\, \text{i.o.} \,\, n \mid\, W^y_{ 1}=x \right) =
\mathbb P \left( X_{\tau^y_n+1}=x \,\,\, \text{i.o.} \,\, n\, \mid X_{\tau^y_1+1}=x \right),
\end{equation*}
to check the recurrence of $(W^y_{n})$ we have to verify that, for all $x \in \mathcal X$,
\begin{equation} \label{recurW}
\mathbb P \left( X_{\tau^y_n+1}=x \,\,\, \text{i.o.} \,\, n\, \mid X_{\tau^y_1+1}=x \right) = 1.
\end{equation}
Fix $x \in \mathcal X$ such that $\mathbb P( X_{\tau^y_1+1}=x ) >0$,
 and choose $\bar x \in \mathcal X$ such that $f_{\bar x}(y) >0$ and $\mathbb P \big ( X_{n+1} = x \mid X_n = \bar x \big) >0$,
(there exists at least one such $\bar x$).
Define the auxiliary sequence of hitting times: 
\begin{eqnarray*}
\sigma^{\bar x,y}_1:=  \inf \{t \geq 0 \mid X_t = \bar x,\, Y_t = y\},\\
\sigma^{\bar x,y}_n := \inf \{t > \sigma^{\bar x,y}_{n-1} \mid X_t = \bar x,\, Y_t = y\}.
\end{eqnarray*}
The hitting times $\big(\sigma^{\bar x,y}_n\big)$ are finite whenever $\big(  \tau^y_n\big)$ are finite,
and the sequence $\big(\sigma^{\bar x,y}_n\big)$
is a subsequence of $\big(  \tau^y_n\big)$, thus
\begin{equation}
( X_{\sigma^{\bar x,y}_n+1}=x) \subseteq  \bigcup_{m\geq n} (X_{\tau^y_m+1}=x),
\end{equation}
and trivially
\begin{equation} \label{EQ_LIMSUP_X}
(\, X_{\sigma^{\bar x,y}_n+1}=x \,\,\, \text{i.o.} \,\, n \,)
\subseteq
(\,X_{\tau_n^y+1}=x \,\,\, \text{i.o.} \,\, n \,).
\end{equation}
The events $\big(X_{\sigma^{\bar x,y}_n +1}= x\big)_n$
are independent under the law
$\mathbb P \big( \cdot \mid X_{\tau^y_1 +1}=x\big)$,
since $\big\{ \big(X_{\sigma^{\bar x,y}_n +1}= x\big)_n,\big(X_{\tau^y_1 +1}=x\big)\big\}$
is a $\mathbb P$-independent set.
In fact, for any $N \in \mathbb N$, and for any choice of $m_1 < \dots<m_N \in \mathbb N$,

\vspace{-3em}
\begin{spacing}{1.5}
\begin{align*}
&\mathbb P \big( X_{\sigma^{\bar x,y}_{m_N}+1}=x, X_{\sigma^{\bar x,y}_{m_{N-1}}+1}=x,\dots,\, X_{\sigma^{\bar x,y}_{m_1}+1}=x ,\,X_{\tau^y_1 +1}=x  \big) &\\
&= \mathbb P \big( X_{\sigma^{\bar x,y}_{m_N}+1}=x,\, X_{\sigma^{\bar x,y}_{m_N}}=\bar x,\,
 \dots,\, X_{\sigma^{\bar x,y}_{m_1}+1}=x,\, X_{\sigma^{\bar x,y}_{m_1}}=\bar x ,\,X_{\tau^y_1 +1}=x    \big) \\
&= \mathbb P \big( X_{\sigma^{\bar x,y}_{m_N}+1}=x\,\mid  X_{\sigma^{\bar x,y}_{m_N}}=\bar x   \big)
\mathbb P \big(  X_{\sigma^{\bar x,y}_{m_N}}=\bar x  \, \mid  X_{\sigma^{\bar x,y}_{m_{N-1}}+1}=x,\,\dots,\, X_{\sigma^{\bar x,y}_{m_1}+1}=x,\, X_{\sigma^{\bar x,y}_{m_1}}=\bar x ,\,X_{\tau^y_1 +1}=x  \big)\\
&\phantom{\qquad \quad} \times\dots\times \mathbb P \big( X_{\sigma^{\bar x,y}_{m_1}+1}=x\,\mid  X_{\sigma^{\bar x,y}_{m_1}}=\bar x   \big)
\times \,\mathbb P \big(X_{\sigma^{\bar x,y}_{m_1}}=\bar x \mid X_{\tau^y_1 +1}=x  \big) \mathbb P \big( X_{\tau^y_1 +1}=x  \big)\\
&= \mathbb P \big( X_{\sigma^{\bar x,y}_{m_N}+1}=x\,  \big)
 \dots \mathbb P \big( X_{\sigma^{\bar x,y}_{m_1}+1}=x\,  \big) \mathbb P \big( X_{\tau^y_1 +1}=x  \big),
\end{align*}
\end{spacing}
\vspace{-1em}
\noindent
where the second equality follows by Remark \ref{REMARK_STRONG_MARKOV_ONE_STOPPING_TIME} in the Appendix,
and the first and last equality follow noting that $X_{\sigma^{\bar x,y}_{m_N}}=\bar x $
for any $\omega \in \Omega$
by definition of $\sigma^{\bar x,y}$.
Note that by definition $\sigma^{\bar x,y}_n > \tau^y_1$ for any $n >1$,
but it could happen $\sigma^{\bar x,y}_1 =\tau^y_1 $,
so the computation above
needs some care for $m_1=1$,
but can be easily recovered also in this case.

The events $\big(X_{\sigma^{\bar x,y}_n +1}= x\big)$ are equiprobable, with strictly positive probability.
By the Borel-Cantelli lemma
\begin{equation}\label{EQ_BOREL_SIGMA}
\mathbb P(\, X_{\sigma^{\bar x,y}_n+1}=x \,\,\, \text{i.o.} \,\, n \, \mid X_{\tau^y_1 +1}=x )=1.
\end{equation}
Equations~(\ref{EQ_LIMSUP_X}) and~(\ref{EQ_BOREL_SIGMA}) taken together give
\begin{align*}
\mathbb P(\,X_{\tau_n^y+1}=x \,\,\, \text{i.o.} \,\, n \, \mid X_{\tau^y_1+1}=x)
\geq   \mathbb P(\, X_{\sigma^{\bar x,y}_n+1}=x \,\,\, \text{i.o.} \,\, n \, \mid X_{\tau^y_1 +1}=x )=1.
\end{align*}
Condition~(\ref{recurW}) is satisfied, thus the recurrence of $(W^y_{n})$ is proved.

\smallskip

\textbf{Step 3} \textit{Verification that the pair $( W_n^y,V_{y,n})$ is a HMM.}
Let us check that the pair $( W_n^y,V_{y,n})$ is as in Definition \ref{DEF_HMM}.
Set
$$    \mathbb P \big(  V_{y,n} = \delta  \mid  W^y_{n} =\varepsilon   \big)=1.$$
For $   \varepsilon \neq x\in \mathcal X$ and $\delta \neq \bar y \in S$,
the pair $(W^y_{n}, V_{y,n})$ inherits the read-out distributions of $(X_n,Y_n)$:
  \begin{equation}\label{READ_OUT_DISTRIBUTION_W_V}
  \mathbb P \big(  V_{y,n} = \bar y   \mid  W^y_{n} =x   \big)=
   \mathbb P \big(  Y_{\tau_n^{y}+1} =\bar  y    \mid  X_{\tau_n^{y}+1} =x   \big) = f_x(\bar y),
  \end{equation}
  see Lemma \ref{LEMMA_READ_OUT_RANDOM_TIMES} in the Appendix.
\noindent
%
Let us verify the
conditional independence property, i.e. that
for any $N \in \mathbb N$
and any $y_1^N \in S^N $ and any $x_1^N \in \mathcal X^N$
\begin{align*}
&\mathbb P \big(  V_{y,1} = y_1, \dots, V_{y,N} = y_N    \mid  W^y_{1} =x_1, \dots, W^y_{N} =x_N   \big) \\
&\phantom{\qquad\qquad} =\prod_{n=1}^N   \mathbb P \big(  V_{y,n} = y_n \mid  W^y_{n} =x_n \big).
\end{align*}
It follows from the direct computation,
\begin{align*}
&\mathbb P \big(  V_{y,1} = y_1 , \dots, V_{y,N} = y_N    \mid  W^y_{1} =x_1, \dots, W^y_{N} =x_N   \big)\\
&=\mathbb P \big(  Y_{\tau_1^{y}+1}  = y_1 , \dots,  Y_{\tau_N^{y}+1}  = y_N
\mid  X_{\tau_1^{y}+1} =x_1, \dots, X_{\tau_N^{y}+1} =x_N   \big)\\
&=\prod_{n=1}^N   \mathbb P \big( Y_{\tau_n^{y}+1}  = y_n \mid  X_{\tau_n^{y}+1} =x_n \big)
=\prod_{n=1}^N   \mathbb P \big(  V_{y,n} = y_n \mid  W^y_{n} =x_n \big),
\end{align*}
where the second equality is a direct consequence of Lemma~\ref{LEMMA_STRONG_INDEPENDENCE_OBSERVATIONS} of the Appendix.
The sequence $(V_{y,n})_{n}$ is therefore a HMM with recurrent underlying Markov chain,
and this concludes the proof of the proposition. $\qed$
\end{proof}


\subsection{Representation of countable mixtures}
\label{Representations of countable Markov mixtures}

In~\cite{DHARMADHIKARI} Dharmadhikari gives a characterization of countable mixtures of i.i.d.
sequences, linking HMMs to the class of exchangeable sequences. The
main result of~\cite{DHARMADHIKARI} can be rephrased as follows
(see Section \ref{tools} for the definitions of exchangeable sequences, mixture of i.i.d. sequences and HMM).
\begin{theorem1} \label{THEOREM_DHARMADHIKARI} (Dharmadhikari)
Let $(Y_n)$ be an exchangeable sequence on a discrete state space $S$.
The sequence
$(Y_n)$ is a \emph{countable} mixture of i.i.d. sequences if and only if $(Y_n)$
is a HMM with recurrent underlying Markov chain.
\end{theorem1}

In the original formulation of
Theorem~\ref{THEOREM_DHARMADHIKARI} the stationarity of the
underlying Markov chain is one of the hypotheses, but close
inspection of the proof in~\cite{DHARMADHIKARI} reveals that only
the absence of transient states is required.

The aim of this
section is to extend the above theorem to partially exchangeable
sequences, i.e. to characterize \emph{countable} mixtures of Markov
chains.
The analog of Theorem
\ref{THEOREM_DHARMADHIKARI} for mixtures of Markov chains is as
follows (we refer to Section \ref{tools} for the definition
of partially exchangeable sequences, mixture of Markov chains and HMMs).
\begin{theorem1}\label{THEOREM_MAIN}
Let $(Y_n)$ be a partially exchangeable sequence on a discrete state space $S$,
with $\mathbb P(Y_0=y_0)=1$
for some $y_0\in S$. The sequence $(Y_n)$ is a \emph{countable} mixture of Markov chains if and only if $(Y_n)$ is a HMM with recurrent underlying Markov chain.
\end{theorem1}
\begin{proof}
\noindent
The standing hypothesis is that $(Y_n)$ is a partially exchangeable sequence.
We first prove that if $(Y_n)$ is a HMM, then it is a countable mixture of Markov chains
i.e., in the notations of
Remark \ref{REMARK_COUNTABLE_STATE_SPACE_MIXTURE_MCS}, $\widetilde P$ takes countably many values.
By the partial exchangeability of $(Y_n)$,
the row $(V_{y, n})$, for any $y \in S$, is exchangeable, and therefore a mixture of i.i.d. sequences.
As proved \emph{e.g.} in Lemma 2.15 of~\cite{ALDOUS} or in Proposition 1.1.4 of~\cite{KALLENBERG_BOOK_2005},
\begin{equation}\label{LLN_exc}
\lim_{N \rightarrow \infty} \frac{1}{N} \sum_{n=1}^N
\mathbb I_{V_{y,n} }(\cdot)= \widetilde p_y(\cdot)\quad \mathbb P-a.s.,
\end{equation}
where the limit has to
be interpreted in the topology of weak convergence, and where $\widetilde p_y$ is the
random probability measure with values in $\mathcal M_{\mathcal S}$ corresponding to
$\widetilde p$ in Definition \ref{DEF_MIXTURES_IID_CONTINUOUS STATE SPACE}.
It follows from the proof of Theorem 1
in~\cite{FORTINI_LADELLI_PETRIS_REGAZZINI}, that the random probability measure
$\widetilde p_y$ in Equation (\ref{LLN_exc}) is the $y$-th row of the random matrix $\widetilde P$.
By Proposition~\ref{LEMMA_SUCCESSIVES MATRIX} above, each row $(V_{y, n})$ is
a HMM, and therefore, by Theorem~\ref{THEOREM_DHARMADHIKARI} above, it is a
countable mixture of i.i.d. sequences. The random probability measure $\widetilde p_y$
is thus concentrated on a countable set, and so is the $y$-th row of $\widetilde P$.
Since this holds for each $y$, the conclusion is that $\widetilde P $ takes countably many values.

To prove the converse one has to show that if
$(Y_n)$ is a given countable mixture of Markov chains,
i.e. if Equation (\ref{EQ_COUNTABLE_MARKOV_MIXTURES_COUNTABLE_STATE_SPACE}) holds
for some countable family $\big(P^h\big)_{h\in H}$, then $(Y_n)$ is a HMM.
We construct a pair $(X_n, \widetilde Y_n)$
satisfying the conditions of Definition \ref{DEF_HMM},
with $(X_n)$
recurrent,
and such that $(\widetilde Y_n )$ and the given $(Y_n)$ have the same distributions. 

The Markov chain $(X_n)$ is defined on the state space\footnote{E.g. ordering the states in first lexical order as follows: $(1,1),(2,1),\dots,(1,2),(2,2),\dots$} $S\times H$, with transition probability  matrix $\mathbf P$, the direct sum of the transition matrices $P^h$,
\begin{eqnarray*}
\mathbf P := \nonumber \left(
\begin{array}{ccccc}
 P^1 & 0 & 0 & \ldots\\
0 &  P^2   & 0 &\ldots \\
\vdots & \vdots & \vdots & \ddots
\end{array} \right),
\end{eqnarray*}
and initial distribution $\pi$ defined, for any $y \in S$ and $h \in H$, as
\begin{equation*}
\pi(y,h)= \left\{
\begin{matrix}
\mu_h& \mbox{ for $y=y_0$ }\\
0  & \mbox{ for $y \neq y_0$,}
\end{matrix} \right.
\end{equation*}
and $\mu_h:= \mathbb P (\widetilde P = P^h)$.
To show that $(X_n) $ is recurrent note that by Theorem 1 in \cite{FORTINI_LADELLI_PETRIS_REGAZZINI},
$(Y_n)$ is conditionally recurrent,
therefore the matrices $\{P^h\}$ in the mixture
correspond to recurrent chains. Since $\mathbf P$ is the
direct sum of such matrices, $(X_n)$ is recurrent.

Consider now a sequence $(\widetilde Y_n)$, with fixed initial state $\widetilde Y_0=y_0$,  conditionally independent given $(X_n)$, and with read-out distributions defined as follows
\begin{equation*}
f_{(x,h)}(y):=\mathbb P \big(\widetilde Y_n=y \mid X_n =(x,h) \big) =\delta_{x,y} ,
\end{equation*}
where $\delta_{\cdot,\cdot}$ is the Kronecker symbol.
Let us compute the finite distributions of $(\widetilde Y_n)$,
for any $N \in \mathbb N$ and any $y_1^N \in (S^*)^{N}$, 

\begin{align*}
&\mathbb P \big(\widetilde  Y_1^N = y_1^N\big)=
\sum_{(x_0^N,\,h_0^N) \in S^{N+1} \times H^{N+1}}  \mathbb P \big(\widetilde  Y_0^N = y_0^N, X_0^N=(x_0^N,h_0^N)  \big)\\
&\phantom{\qquad}=\sum_{(x_0^N,\,h_0^N) \in S^{N+1} \times H^{N+1}} \mathbb P (X_0 = (x_0,h_0))\, \prod_{n=0}^N   \mathbb P \big(\widetilde  Y_n = y_n \mid X_n=(x_n,h_n) \big) \\
&\phantom{\qquad\qquad\qquad}\times \prod_{n=0}^{N-1} \mathbb P \big(  X_{n+1}=(x_{n+1},h_{n+1})\mid X_n=(x_n,h_n) \big) \\
&\phantom{\qquad}=  \sum_{(x_0^N,\,h_0^N) \in S^{N+1} \times H^{N+1}}
\pi(x_0,h_0) \prod_{n=0}^N     \,f_{(x_n,h_n)}(y_n) \,   \prod_{n=0}^{N-1} \mathbf P_{(x_{n},h_{n})(x_{n+1},h_{n+1})}\\
&\phantom{\qquad}= \,\, \sum_{h \in H} \,\, \mu_h P^{h}_{y_0 y_{1}} \dots P^{h}_{y_{N-1} y_{N}} ,
\end{align*}
where the second equality follows from the conditional independence of $(\widetilde Y_n)$ given $(X_n)$, and the fourth from the definition of the read-out densities and by the block structure of $\mathbf P$.
Comparing~(\ref{EQ_COUNTABLE_MARKOV_MIXTURES_COUNTABLE_STATE_SPACE}) with the last
expression, we have that $(Y_n)$ and $(\widetilde Y_n)$ have the same distributions,
thus $(Y_n)$ is a HMM with recurrent underlying Markov chain and the theorem is proved.
\end{proof}

 Note that, if $S$ finite, in Theorems \ref{THEOREM_DHARMADHIKARI} and \ref{THEOREM_MAIN} the state space of the underlying Markov chain is \emph{finite} if and only if the mixture is
\emph{finite}.

\section{Countable Markov mixtures with Polish state space}\label{Polish state spaces}
\setcounter{equation}{0}
In this section $S$ is a Polish space.
\subsection{The successors array}

The proposition below is the analog of Proposition \ref{LEMMA_SUCCESSIVES MATRIX}
for uncountable state space $S$ (the definitions of HMM and of successors array are in Section \ref{tools}).
\begin{proposition1} \label{PROPOSITION_SUCCESSIVES MATRIX_CONTINUOUS_STATE_SPACE}
Let $(Y_n)$ be a HMM on a Polish space $S$ with recurrent underlying Markov
chain, then each row of the successors array $(V_{j,n})$ is a HMM with recurrent
underlying Markov chain.
\end{proposition1}

\begin{proof} 
Let $(X_n)$ be the underlying Markov chain of $(Y_n)$.
Consider the partition $\mathcal E = (E_j)_{j \geq 0}$ of $S^*$, and for any element
$E_j$ of the partition define
\begin{equation*}
\tau^{E_j}_1:=\inf\{t\geq 0 \mid Y_t \in E_j\}, \quad \tau^{E_j}_n:=\inf\{t> \tau^{E_j}_{n-1} \mid Y_t \in E_j\}.
\end{equation*}
The proof can be carried out exactly as
the proof of Proposition \ref{LEMMA_SUCCESSIVES MATRIX}, substituting
$\tau^{y}_n$ there with $\tau^{E_j}_n$, and $\sigma^{\bar x,y}_n$ with $\sigma^{\bar x,E_j}_n$, defined below
\begin{align*}
\sigma^{\bar x,E_j}_1 &:=  \inf \{t \geq 0 \mid X_t = \bar x,\, Y_t \in E_j\},\\
\sigma^{\bar x,E_j}_n &:= \inf \{t > \sigma^{\bar x,y}_{n-1} \mid X_t = \bar x,\, Y_t  \in E_j\},
\end{align*}
where $\bar x \in \mathcal X$ is such that $f_{\bar x}(E_j) >0$ and $\mathbb P \big ( X_{n+1} = x \mid X_n = \bar x \big) >0$,
($x $ has the same role as in Equation (\ref{recurW})).
\end{proof}

\subsection{Representation of countable mixtures}

This subsection mirrors Subsection \ref{Representations of countable Markov mixtures} for the case of Polish state space $S$.
Theorem \ref{THEOREM_DHARMADHIKARI_CONTINUOUS_STATE_SPACE} below extends Theorem~\ref{THEOREM_DHARMADHIKARI} to Polish state spaces. To the best of our knowledge the extension is not available in the literature. Based on Theorem \ref{THEOREM_DHARMADHIKARI_CONTINUOUS_STATE_SPACE} we prove Theorem \ref{THEOREM_MAIN_CONTINUOUS_STATE_SPACE} which is the counterpart of Theorem \ref{THEOREM_MAIN} and the main result of the subsection.

Note that Theorem~\ref{THEOREM_DHARMADHIKARI}, i.e. Dharmadhikari's original result \cite{DHARMADHIKARI},  can not be directly generalized as it relies on a definition of HMMs unsuitable for general state spaces.

\subsubsection{Representation of countable i.i.d. mixtures}
Exchangeable sequences, mixture of i.i.d. sequences, and HMMs are defined in Section \ref{tools}.
\begin{theorem1} \label{THEOREM_DHARMADHIKARI_CONTINUOUS_STATE_SPACE}
Let $(Y_n)$ be an exchangeable sequence on a Polish space. The sequence $(Y_n)$ is a
\emph{countable} mixture of i.i.d. sequences if and only if $(Y_n)$
is a HMM with recurrent underlying Markov chain.
\end{theorem1}

\begin{proof}
If $(Y_n)$ is a countable mixture of i.i.d.,
then it is a HMM by Remark \ref{RMK_HMM_MIXT_IID}. To prove the converse let $(Y_n )$ be an exchangeable HMM, whose recurrent underlying Markov chain $(X_n)$ has transition probability matrix $P$ and initial distribution $\pi$. The Markov chain $(X_n)$ has no transient states, but possibly more than one recurrence class.
As noted in \cite{DHARMADHIKARI}, by the exchangeability of $(Y_n )$,
one can substitute $P$ with the Ces\`aro limit
$P^*:= \lim_{n \rightarrow \infty } 1/n \sum_{k=1}^n P^k$,
where $P^k$ is the $k$-power of $P$.
By the ergodic theorem $P^*$ has a block structure, being the direct sum of matrices $P_h$ with
identical rows, one block $P_h$ for each recurrence class.
By Lemma \ref{LEMMA_LUMPING_CONTINUOUS_READ_OUT}
of the Appendix $(Y_n)$ is a countable mixture of i.i.d. sequences.
\end{proof}

\subsubsection{Representation of countable Markov mixtures}

See Section \ref{tools} for the definitions of partially exchangeable sequence,
mixture of Markov chains and HMM.
\begin{theorem1}\label{THEOREM_MAIN_CONTINUOUS_STATE_SPACE}
Let $(Y_n)$ be a partially exchangeable sequence on a Polish space
with $\mathbb P(Y_0=y_0)=1$, for some $y_0 \in E_1$.
The sequence $(Y_n)$ is a \textit{countable} mixture of homogeneous Markov chains if and only if $(Y_n)$ is a HMM with recurrent underlying Markov chain.
\end{theorem1}
\begin{proof}
The partially exchangeability of $(Y_n)$
is a standing hypothesis.
Let $(Y_n)$ be a HMM with recurrent underlying Markov chain.
To prove that $(Y_n)$ is a countable mixture of Markov chains imitate
the proof of Theorem \ref{THEOREM_MAIN}.
Note first that $\frac{1}{N} \sum_{n=1}^N\mathbb I_{V_{j,n}} (\cdot)\rightarrow \theta_j(\cdot)$,
where $ \theta_j$ is a probability measure on $\mathcal S$.
As in the proof of Theorem 4 in \cite{FORTINI_LADELLI_PETRIS_REGAZZINI},
for any $j \in \mathbb N$ define $ \widetilde t (j,\cdot):=  \theta_j (\cdot) $.
To conclude use Proposition \ref{PROPOSITION_SUCCESSIVES MATRIX_CONTINUOUS_STATE_SPACE} and Theorem \ref{THEOREM_DHARMADHIKARI_CONTINUOUS_STATE_SPACE}.

For the converse assume that $(Y_n)$ is a countable mixture of Markov chains, with random kernel $\widetilde t$ taking values
$\big(t_h\big)_{h \in H}$ and with $\mu_h= \mathbb P (\widetilde t = t_h)$,
and finite distributions as in Equation (\ref{EQ_COUNTABLE_MARKOV_MIXTURES_EXPLICIT}).
To prove that $(Y_n)$ is a HMM with recurrent underlying Markov chain we construct a recurrent Markov chain $(X_n)$
and a sequence $(\widetilde Y_n)$ satisfying the first two conditions in Definition \ref{DEF_HMM},
then showing that $(\widetilde Y_n)$ has the same distributions of $(Y_n)$
\footnote{
The construction used for the proof of Theorem~\ref{THEOREM_MAIN} can not be used here,
in fact the Markov chain $(X_n)$ there
takes values in the
product space $H\times S$,
 which can be now uncountable, while we need a discrete underlying Markov chain.}.
Consider thus a Markov chain $(X_n)$ taking values in $H \times \mathbb N_0 \times \mathbb N_0 $, with components
$X_n=(h_n,i_{n},j_{n})$ representing the index of the running chain in the mixture,
the discretized value of $Y_{n-1}$ (i.e. the elements of the partition to which $Y_{n-1}$ belongs),
 and the discretized value of $Y_n$ respectively.
The initial distribution of $(X_n)$ is taken to be
\[
\mathbb P \big(  X_0 = (h, i,j) \big) \,=\,
 \mu_h  \frac{1}{2^{i+1}}  \delta_{j,1}\, ,
\]
where $\delta_{\cdot,\cdot}$ is again the Kronecker symbol,
and its transition probabilities
$$
\mathbb P \Big(X_{n} = (h_n,i_{n},j_{n}) \mid X_{n-1} = (h_{n-1},i_{n-1},j_{n-1}) \Big) \,=\, \delta_{h_{n-1},h_n }\delta_{i_{n}, j_{n-1} } t_{h_n}(i_{n}, E_{j_{n}}).
$$

The Markov chain $(X_n)$ is recurrent since
the kernels $t_h$ correspond to recurrent Markov chains
by Theorem 4 in \cite{FORTINI_LADELLI_PETRIS_REGAZZINI}.
Consider now a sequence $(\widetilde Y_n)$ jointly distributed with $(X_n)$, with fixed initial value $\widetilde Y_0 =y_0$, conditionally independent given $(X_n)$, and with read-out distributions defined as follows for any $A \in \mathcal S$
\begin{align*}
\mathbb P \big(  \widetilde Y_n\in A    \mid   X_n = (h,i,j)     \big) & =
\left\{
\begin{array}{cc}
 0 & \mbox{ for $t_h(i, E_{j})=0 $}  \\
 \frac{1}{t_h(i, E_{j})}\int_{A \cap E_j}t_h(i,d y)  & \mbox{ for $t_h(i, E_{j})\neq 0$}.
\end{array}
\right.
\end{align*}

For any $N \in \mathbb N$,
and any $A_1, \dots, A_N \in \mathcal S^*$,
the distributions of $(\widetilde Y_n)$ are computed as follows
\begin{align*}
&\mathbb P \big( \widetilde  Y_1 \in A_1 , \dots, \widetilde Y_N \in A_N  \big)\\
&= \sum_{h_1^N \in H^N, \,i_1^N \in (\mathbb N_0)^{N},\, j_1^N \in (\mathbb N_0)^{N}}
\mathbb P \big(  \widetilde Y_1 \in A_1, \dots,\widetilde Y_N \in A_N   , \,    X_1 = (h_1,i_1,j_1), \dots, X_N =(h_N, i_N,j_N) \big)\\
&= \sum_{h_1^N \in H^N,\, i_1^N \in (\mathbb N_0)^{N},\, j_1^N \in (\mathbb N_0)^{N}}
\mathbb P \big(  \widetilde Y_1 \in A_1, \dots,\widetilde Y_N \in A_N  \mid   X_1 = (h_1,i_1,j_1), \dots, X_N =(h_N,i_N,j_N) \big)   \\
&\phantom{\qquad\qquad}\times \mathbb P \big( X_1 = (h_1,i_1,j_1), \dots, X_N =(h_N,i_N,j_N) \big) \\
&= \sum_{h_0^N \in H^N,\, i_0^N\in (\mathbb N_0)^{N+1}, \,j_0^N \in (\mathbb N_0)^{N+1}}
\prod_{n=1}^N  \mathbb P \big(  \widetilde Y_n \in A_n \mid   X_n = (h_n,i_n,j_n)\big)   \\
&\phantom{\qquad\qquad}\times\prod_{n=1}^{N}  \mathbb P \big( X_{n} = (h_n,i_{n},j_{n}) \mid X_{n-1} =(h_{n-1},i_{n-1},j_{n-1}) \big)
\,\mathbb P \big( X_0 = (h_0,i_{0},j_{0}) \big) \\
&=\sum_{h_0^N \in H^N,\, i_0^N \in (\mathbb N_0)^{N+1},\, j_0^N\in (\mathbb N_0)^{N+1}} \prod_{n=1}^N  \frac{1}{t_{h_n}(i_n, E_{j_n})}\int_{A_n \cap E_{j_n}}t_{h_n}(i_n,d y) \\
&\phantom{\qquad\qquad}\times \prod_{n=1}^{N}  \delta_{h_{n-1}, h_{n} }\, \delta_{i_{n}, j_{n-1} } t_{h_n}(i_{n}, E_{j_{n}})
 \,\mu_{h_0} \delta_{j_0,1}   \frac{1}{2^{i_0+1}}\\
&=\sum_{h \in H,\, i_0 \in \mathbb N_0,\, j_0^N \in (\mathbb N_0)^{N+1}} \prod_{n=1}^N  \int_{A_n \cap E_{j_n}}t_{h}(j_{n-1},d y)
 \,\mu_{h} \delta_{j_0,1}   \frac{1}{2^{i_0+1}}\\
&=\sum_{h \in H,  \,j_1^N \in (\mathbb N_0)^{N}} \mu_{h}\int_{A_1 \cap E_{j_1}}t_{h}(1,d y)  \prod_{n=2}^N  \int_{A_n \cap E_{j_n}}t_{h}(j_{n-1},d y).
\end{align*}
 Comparing the expression above with Equation (\ref{EQ_COUNTABLE_MARKOV_MIXTURES_EXPLICIT}),
one concludes that the distributions of $(\widetilde Y_n)$ coincide with those of $(Y_n)$, therefore proving that $(Y_n)$ is a HMM.
\end{proof}

\section{Concluding remarks} \label{sec:conclus}
\setcounter{equation}{0}
Throughout the paper we referred to the notion of partial
exchangeability originally given by de Finetti and to the
corresponding representation theorem as given
in~\cite{FORTINI_LADELLI_PETRIS_REGAZZINI}.
For discrete state space
partial exchangeability can be defined in a slightly different way,
and a representation theorem in this alternative framework 
is proved in \cite{DIACONIS_FREEDMAN_1980}.
According to \cite{DIACONIS_FREEDMAN_1980},
a sequence of random variables is
partially exchangeable if the probability is invariant under all
permutations of a string that preserves the first value and the
transition counts between any couple of states.
A characterization of countable mixtures
of Markov chains can be given also in the setup of~\cite{DIACONIS_FREEDMAN_1980},
using different mathematical tools. 
The result is in~\cite{FINESSO_PROSDOCIMI_ECC},
but for a complete proof see~\cite{PROSDOCIMI}.
By the same token the characterization of countable
mixtures of Markov chains of order $k$ holds true, for the proof
see~\cite{PROSDOCIMI}. Unfortunately the approach of~\cite{FINESSO_PROSDOCIMI_ECC}
and~\cite{PROSDOCIMI} does not readily generalize to Polish state space.

Based on the results in~\cite{FORTINI_LADELLI_PETRIS_REGAZZINI},
a de Finetti's type representation theorem for mixtures of semi-Markov
processes have been proved in~\cite{EPIFANI_FORTINI_LADELLI}.
The authors are confident that a characterization of countable mixtures of semi-Markov processes
in terms of HMMs
can be given properly adapting the proof of
Proposition~\ref{PROPOSITION_SUCCESSIVES MATRIX_CONTINUOUS_STATE_SPACE}
and Theorem~\ref{THEOREM_MAIN_CONTINUOUS_STATE_SPACE}.

\section{Appendix}
\setcounter{equation}{0}
\subsection{Strong Markov and strong conditional independence for HMMs}
This section contains some useful properties of HMMs.

\begin{lemma1}\label{LEMMA_HMM_SPLITTING}(Splitting property)
Let $(Y_n)$ be a HMM with underlying Markov chain $(X_n)$.
Then the pair $(X_n,Y_n)$ is a Markov chain. Moreover
for any $N \in \mathbb N$,
for any $x_1, \dots, x_N \in \mathcal X$, and $S_1, \dots, S_N \in \mathcal S $
such that
$\mathbb P \big(  X_1^{N-1}= x_1^{N-1} , \,Y_1^{N-1}  \in S_1^{N-1} \big)>0$
we have
\begin{align*}
&\mathbb P \big(  X_N = x, \,  Y_N \in S_N \mid X_1^{N-1}= x_1^{N-1} , \,Y_1^{N-1}  \in S_1^{N-1} \big)&\\
&\phantom{\qquad}= \mathbb P \big(   X_N = x, \,  Y_N \in S_N  \mid X_{N-1}= x_{N-1}\big).
\end{align*}
\end{lemma1}
\begin{proof}
\begin{spacing}{1.5}\begin{align*}
&\mathbb P \big(  X_N = x, \,  Y_N \in S_N \mid X_1^{N-1}= x_1^{N-1} , \,Y_1^{N-1}  \in S_1^{N-1} \big)\\
&\phantom{\quad}= \frac{\mathbb P \big(  X_1^{N}= x_1^{N} , \,Y_1^{N}  \in S_1^{N} \big)}
{\mathbb P \big(    X_1^{N-1}= x_1^{N-1} , \,Y_1^{N-1}  \in S_1^{N-1} \big)}\\
&\phantom{\quad}= \frac{\prod_{n=1}^N  \mathbb P \big(  Y_n \in S_n\mid X_n = x_n \big) \,\, \mathbb P \big( X_1^{N} = x_1^{N}  \big) }
{ \prod_{n=1}^{N-1}  \mathbb P \big(  Y_n \in S_n \mid X_n = x_n  \big)
\,\,   \mathbb P \big( X_1^{N-1} = x_1^{N-1}  \big)  }\\
&\phantom{\quad}= \mathbb P \big( Y_N \in S_N \mid X_N=x_N\big) \,\, \mathbb P \big( X_N = x_N \mid X_{N-1} = x_{N-1} \big) \\
&\phantom{\quad}=  \mathbb P \big( Y_N \in S_N \mid X_N=x_N,\, X_{N-1} = x_{N-1}\big) \,\, \mathbb P \big( X_N = x_N \mid X_{N-1} = x_{N-1} \big) \\
&\phantom{\quad}= \mathbb P \big( Y_N \in S_N,\, X_N=x_N \mid X_{N-1} = x_{N-1}\big),
\end{align*}
\end{spacing}
\vspace{-0.5em}
\noindent
where the second and fourth equality follow by the conditional independence of the observations $(Y_n)$ in the definition of HMM.
\end{proof}

\begin{lemma1}\label{LEMMA_STRONG_MARKOV_ONE_STOPPING_TIME}
(Strong splitting property)
Let $(Y_n)$ be a HMM with underlying Markov chain $(X_n)$, and
$\gamma$ be a stopping time for $(X_n,Y_n)$,
then, for any $x, \widetilde x, \bar x \in \mathcal X,$ and any $ S_1,S_2,S_3 \in \mathcal S$
such that
 $\mathbb P \big(  X_{\gamma}= \widetilde x , \,Y_{\gamma } \in S_2,\, X_{\gamma \wedge n} = \bar x,\, Y_{\gamma \wedge n} \in S_1 \big) >0$
it holds that
\begin{align}\label{EQ_STRONG_MARKOV_ONE_STOPPING_TIME}
&\mathbb P \big(  X_{\gamma+k} = x, \,  Y_{\gamma+k}\in S_3 \mid X_{\gamma}= \widetilde x , \,Y_{\gamma } \in S_2,\, X_{\gamma \wedge n} = \bar x,\, Y_{\gamma \wedge n} \in S_1 \big) \nonumber \\
&\phantom{\qquad}= \mathbb P \big(   X_{\gamma+k} = x, \,  Y_{\gamma+k} \in S_3 \mid X_{\gamma}= \widetilde x \big).
\end{align}
\end{lemma1}
\begin{proof}
We manipulate separately the left-hand side (LHS)
 and the right-hand side (RHS) of Equation~(\ref{EQ_STRONG_MARKOV_ONE_STOPPING_TIME}). For readability
denote $C_r:=\big( \gamma=r,\,  X_{r}= \widetilde x , \,Y_{r }\in S_2,\, X_{r \wedge n} = \bar x,\,
Y_{r \wedge n} \in S_1 \big)$.
Applying Lemma \ref{LEMMA_HMM_SPLITTING}, the numerator of the conditional probability on the LHS of Equation~(\ref{EQ_STRONG_MARKOV_ONE_STOPPING_TIME}) is
\vspace{-2em}\begin{spacing}{1.5}
\begin{align*}
&\mathbb P \big(  X_{\gamma+k} = x, \,  Y_{\gamma+k} \in S_3 , \,X_{\gamma}= \widetilde x , \,Y_{\gamma}\in S_2,\, X_{\gamma \wedge n} = \bar x,\, Y_{\gamma \wedge n} \in S_1 \big)\\
&=\textstyle{\sum_{r\geq 1}}\mathbb P \big(  X_{r+k} = x, \,  Y_{r+k} \in S_3 \mid \,C_r \big)
 \mathbb P(C_r)\\
&=\textstyle{\sum_{r\geq 1}}\mathbb P \big(  X_{r+k} = x, \,  Y_{r+k} \in S_3 \mid \,  X_{r}= \widetilde x   \big) \mathbb P(C_r)\\
&=  \textstyle{\sum_{r\geq 1}}    \mathbb P \big(    Y_{r+k} \in S_3 \mid \,  X_{r+k} = x \big)
  \mathbb P \big(  X_{r+k} = x \,\mid  X_{r}= \widetilde x   \big) \mathbb P(C_r)\\
&= f_x(S_3)  P_{\widetilde x, x}^{(k)}  \,\,\textstyle{\sum_{r\geq 1}} \mathbb P(C_r)\\
%
&=f_x( S_3)  P_{\widetilde x, x}^{(k)}  \,\, \mathbb P \big( X_{\gamma}= \widetilde x , \,Y_{\gamma }\in S_2,\, X_{\gamma \wedge n} = \bar x,\, Y_{\gamma \wedge n} \in S_1 \big)   ,
\end{align*}
\end{spacing}
\noindent
where $ P_{\widetilde x, x}^{(k)}$ is the ${\widetilde x, x}$-entry of
the $k$-step transition matrix of the Markov chain $(X_n)$.
The numerator of the conditional probability on the RHS of Equation (\ref{EQ_STRONG_MARKOV_ONE_STOPPING_TIME}), again applying Lemma \ref{LEMMA_HMM_SPLITTING}, is
\begin{spacing}{1.5}
\begin{align*}
&\mathbb P \big(   X_{\gamma+k} = x, \,  Y_{\gamma+k} \in S_3,\,  X_{\gamma}= \widetilde x \big)\\
%
&=    \textstyle{\sum_{r\geq 1}}   \mathbb P \big(  X_{r+k} = x, \,  Y_{r+k} \in S_3 \, \mid\,\gamma=r,\,  X_{r}= \widetilde x \big)\, \mathbb P \big(  \gamma=r,\,  X_{r}= \widetilde x \big)\\
&=  \textstyle{\sum_{r\geq 1}}  \mathbb P \big(  X_{r+k} = x, \,  Y_{r+k} \in S_3\, \mid\, X_{r}= \widetilde x \big)\, \mathbb P \big(  \gamma=r,\,  X_{r} = \widetilde x \big) \\
&=f_x(S_3)  P_{\widetilde x, x}^{(k)}  \mathbb P \big(  X_{\gamma}= \widetilde x \big).
\end{align*}
\end{spacing}
\noindent
The lemma is proved comparing the expressions of the LHS and the RHS derived above.
\end{proof}

Taking $S_1= S_2= S_3=S$ we have
\begin{remark1}\label{REMARK_STRONG_MARKOV_ONE_STOPPING_TIME}
Let $(X_n), \gamma$ be as in Lemma \ref{LEMMA_STRONG_MARKOV_ONE_STOPPING_TIME}
then, for any $x, \widetilde x, \bar x \in \mathcal X,$
 such that $\mathbb P \big(  X_{\gamma}= \widetilde x , \, X_{\gamma \wedge n} = \bar x \big) >0$
it holds that
\begin{align}%
&\mathbb P \big(  X_{\gamma+k} = x, \, \mid X_{\gamma}= \widetilde x , \, X_{\gamma \wedge n} = \bar x\,
 \big)= \mathbb P \big(   X_{\gamma+k} = x, \,  \mid X_{\gamma}= \widetilde x \big).
\end{align}
\end{remark1}

\begin{definition1}\label{DEF_HITTING_TIMES}
Let $(Y_n)$ be a HMM with underlying Markov chain $(X_n)$,
and let $A \subset \mathcal X \times \mathcal S$.
We say that the sequence of random times $\big(\gamma_n\big)_{n \geq 1}$ is a
sequence of \emph{hitting times of $A$} if
\begin{align*}
\gamma_1 &:= \inf\{ t \geq 0 \mid  (X_t,Y_t) \in A \},\\
\gamma_n &:= \inf\{t > \gamma_{n-1} \mid    (X_t,Y_t) \in A \}.
\end{align*}
\end{definition1}
\begin{lemma1}\label{LEMMA_n_STOPPING_TIMES}  {(Generalized strong splitting property)}
Let $(Y_n)$ be a HMM with underlying Markov chain $(X_n)$.
Let $(\gamma_n)$ be a sequence of hitting times of $A$ for $(X_n,Y_n)$,
where $A \subset \mathcal X \times \mathcal S$.
Then for any $N$, and any $(x_1,S_1), \dots,(x_N, S_N) \in A $
such that
$\mathbb P\big(X^{\gamma_{N-1}}_{\gamma_1} = x^{{N-1}}_1, \,Y^{\gamma_{N-1}}_{\gamma_1} \in S^{{N-1}}_1 \big)>0
$
it holds
\begin{align*}\label{EQ_TWO_STOPPING_TIMES}
&\mathbb P\big(X_{\gamma_N} = x_{N}, \, Y_{\gamma_N} \in S_{N}
\mid X^{\gamma_{N-1}}_{\gamma_1} = x^{{N-1}}_1, \,Y^{\gamma_{N-1}}_{\gamma_1} \in S^{{N-1}}_1 \big)
\\
&\phantom{\qquad\qquad}= \mathbb P \big(   X_{\gamma_N} = x_{N}, \, Y_{\gamma_N} \in S_{N}  \mid X_{\gamma_{N-1}}=  x_{{N-1}}\big).
\end{align*}
\end{lemma1}
\begin{proof}
Denote with $A^c$ the complement of $A$ in $\mathcal X \times \mathcal S$,
and with $(A^c)^r$ the $r$-th fold Cartesian product of $A^c$.
Let $B := \left( X^{\gamma_{N-1}}_{\gamma_1}=  x^{{N-1}}_1, \, Y^{\gamma_{N-1}}_{\gamma_1} \in S^{{N-1}}_1 \right)$.
Applying Lemma \ref{LEMMA_STRONG_MARKOV_ONE_STOPPING_TIME} in the third equality below,
the numerator of the conditional probability on the LHS is
\vspace{-2em}
\begin{spacing}{1.5}
\begin{align*}
\left. \right.&\mathbb P \big( X^{\gamma_N}_{\gamma_1} = x^N_1, \,Y^{\gamma_N}_{\gamma_1} \in S^N_1 \big) \\
&= \textstyle{\sum_{r \geq 1}}
\mathbb P \big(\gamma_{N}= \gamma_{N-1}+r, \, X_{\gamma_{N-1}+r}=  x_{{N}}  , \,Y_{\gamma_{N-1}+r} \in S_{{N}}\, \mid B \big)
 \mathbb P \big( B \big)\\
&= \textstyle{\sum_{r \geq 1}}
\mathbb P \Big( X_{\gamma_{N-1}+r}=x_{{N}}, \,Y_{\gamma_{N-1}+r} \in S_{{N}}, \,
\big(X^{\gamma_{N-1}+r-1}_{\gamma_{N-1}+1}, \,Y^{\gamma_{N-1}+r-1}_{\gamma_{N-1}+1} \big) \in (A^c)^r \mid  B \Big) \mathbb P \big(  B \big)\,  \\
%
%
&= \textstyle{\sum_{r \geq 1}}
\mathbb P \Big( X_{\gamma_{N-1}+r}=x_{{N}}, \,Y_{\gamma_{N-1}+r} \in S_{{N}},
\big(X^{\gamma_{N-1}+r-1}_{\gamma_{N-1}+1}, \,Y^{\gamma_{N-1}+r-1}_{\gamma_{N-1}+1} \big) \in (A^c)^r\, \mid  X_{\gamma_{N-1}}\!\!=\!  x_{{N-1}}  \Big) \mathbb P \big(  B \big) \\
&= \textstyle{\sum_{r \geq 1}}
\mathbb P \big(\gamma_{N}= \gamma_{N-1}+r,\, X_{\gamma_{N-1}+r}=  x_{{N}}  , \,Y_{\gamma_{N-1}+r} \in S_{{N}}\,   \mid X_{\gamma_{N-1}}=  x_{{N-1}} \big) \mathbb P \big(  B\big)\\
&= \mathbb P \big( X_{\gamma_{N}}=  x_{{N}}  , \,Y_{\gamma_{N}} \in S_{{N}}\, \mid X_{\gamma_{N-1}}=  x_{{N-1}}   \big)
 \mathbb P \big(B\big),
\end{align*}
\end{spacing}

\vspace{-1em}\noindent
and dividing by $\mathbb P(B)$ the lemma is proved.
\end{proof}

\begin{remark1}\label{REMARK_STRONG_SPLITTING}
By the same token, for any $(x_1,S_1), \dots,(x_N, S_N) \in \mathcal X \times \mathcal S $,
\begin{align*}
\left. \right. &\mathbb P \big(  X_{\gamma_N+1} = x_{N}, \, Y_{\gamma_N+1} \in S_{N}
\mid  X^{\gamma_{N-1}+1}_{\gamma_1+1}=  x^{{N-1}}_1, \, Y^{\gamma_{N-1}+1}_{\gamma_1+1} \in S^{{N-1}}_1 \big) \\
&\phantom{\qquad}=\mathbb P \big(   X_{\gamma_N +1} = x_{N}, \, Y_{\gamma_N +1} \in S_{N}  \mid X_{\gamma_{N-1} +1}=  x_{{N-1}}\big).
\end{align*}
\end{remark1}

\noindent
Taking $S_1= \dots=S_N =S$, Remark~\ref{REMARK_STRONG_SPLITTING} gives
\begin{remark1}\label{REMARK_STRONG_STRONG_MARKOV}
For any $x_1, \dots, x_N \in \mathcal X$,
\begin{align*}
\mathbb P \big(  X_{\gamma_N+1} = x_{N} \mid X^{\gamma_{N-1}+1}_{\gamma_1+1}=  x^{{N-1}}_1 \big)
= \mathbb P \big(   X_{\gamma_N+1} = x_{N} \mid X_{\gamma_{N-1}+1}=  x_{{N-1}}\big).
\end{align*}
\end{remark1}

As a consequence of the conditional independence property of HMMs we have
\begin{lemma1}\label{LEMMA_READ_OUT_RANDOM_TIMES}
Let $(Y_n)$ be a HMM with underlying Markov chain $(X_n)$.
Then for any $N \in \mathbb N$,
for any $x_1, \dots, x_N \in \mathcal X$, and $S_1, \dots, S_N \in \mathcal S $
such that
$\mathbb P \big(  X_1^{N}= x_1^{N} , \,Y_1^{N-1}  \in S_1^{N-1} \big)>0$
we have
\begin{equation}\label{EQ_ELABORATION_STRONG_SPLITTING}
\mathbb P \left( Y_N \in  S_N \, \mid  X_1^{N}= x_1^{N}  , \,Y_1^{N-1}  \in S_1^{N-1}\right)
=
\mathbb P \left( Y_N \in  S_N \, \mid  X_N= x_N  \right).
\end{equation}
Moreover let $\sigma, \tau$ two stopping times for $(X_n,Y_n)$
such that $\sigma < \tau$.
Then for any $\bar S \in \mathcal S^*$, and any $x_1,x_2 \in \mathcal X$ we have
 \begin{equation}\label{EQ_READ_OUT_STOPPING_TIME}
\mathbb P \left( Y_{\tau+1} \in  \bar S \, \mid  X_{\tau+1}= x_2 \right)
=
f_{x_2}(\bar S),
\end{equation}
 \begin{equation}\label{EQ_STRONG_READ_OUT}
\mathbb P \left( Y_{\tau+1} \in  \bar S \, \mid  X_{\tau+1}= x_2  , \, X_{\sigma+1}= x_1 \right)
=
\mathbb P \left( Y_{\tau+1} \in  \bar S \, \mid  X_{\tau+1}= x_2   \right).
\end{equation}
\end{lemma1}

%
%
%
%
%
\begin{proof}
Equation (\ref{EQ_ELABORATION_STRONG_SPLITTING}) can be easily proved
using the conditional independence property.
Equation (\ref{EQ_READ_OUT_STOPPING_TIME}) can be seen disintegrating the stopping time $\tau+1$,
\begin{align*}
&\mathbb P \big(  Y_{\tau+1} \in \bar  S ,\,    X_{\tau+1} =x   \big) \\
&=  \sum_{m \in \mathbb N} \mathbb P \big( \tau+1=m,\, Y_{m} \in \bar  S ,\,  X_{m} =x   \big)\\
&= \sum_{m \in \mathbb N}  \mathbb P \big(  Y_{m} \in \bar  S \mid \tau+1=m,\, X_{m} =x   \big)
  \mathbb P \big(  \tau+1=m,\,X_{m} =x   \big)\\
&= \sum_{m \in \mathbb N}  \mathbb P \big(Y_{m} \in \bar  S \mid  X_{m} =x   \big)
  \mathbb P \big( \tau+1=m,\, X_{m} =x   \big)\\
&= f_x(\bar S) \sum_{m \in \mathbb N} \mathbb P \big( \tau+1=m,\, X_{m} =x   \big)
= f_x(\bar S) \mathbb P \big(  X_{\tau+1} =x   \big),
\end{align*}
where the third equality follows by Equation (\ref{EQ_ELABORATION_STRONG_SPLITTING}),
and the result then follows by definition of conditional probability.

\noindent
To get Equation (\ref{EQ_STRONG_READ_OUT}) statement write
\begin{align*}
 & \mathbb P \left(    Y_{\tau+1} \in  \bar S,  \, X_{\tau+1}= x_2  , \, X_{\sigma+1}= x_1 \right)\\
&= \sum_{m >n} \sum_{n \in \mathbb N}
\mathbb P \left( \tau =m, \, \sigma=n, \,Y_{m+1} \in  \bar S, \, X_{m+1}= x_2  , \, X_{n+1}= x_1 \right)\\
&= \sum_{m >n} \sum_{n \in \mathbb N}
\mathbb P \left( Y_{m+1} \in  \bar S \mid \tau =m, \, \sigma=n,\, X_{m+1}= x_2  , \, X_{n+1}= x_1 \right)
\mathbb P \left( \tau =m, \, \sigma=n,\, X_{m+1}= x_2  , \, X_{n+1}= x_1 \right)\\
&=
\sum_{m >n} \sum_{n \in \mathbb N}
\mathbb P \left( Y_{m+1} \in  \bar S\, \mid\, X_{m+1}= x_2   \right)
\mathbb P \left( \tau =m, \, \sigma=n,\, X_{m+1}= x_2  , \, X_{n+1}= x_1 \right)\\
&= f_{x_2} (\bar S)  \sum_{m >n} \sum_{n \in \mathbb N}
\mathbb P \left( \tau =m, \, \sigma=n,\, X_{m+1}= x_2  , \, X_{n+1}= x_1 \right)\\
&
=f_{x_2} (\bar S)
\mathbb P \left(  X_{\tau+1}= x_2  , \, X_{\sigma+1}= x_1 \right),
\end{align*}
where the third equality follows  by Equation (\ref{EQ_ELABORATION_STRONG_SPLITTING}),
and the result follows by the first statement and again by definition
of conditional probability.
\end{proof}

\begin{lemma1}\label{LEMMA_STRONG_INDEPENDENCE_OBSERVATIONS}
Let $(Y_n)$ be a HMM with underlying Markov chain $(X_n)$, and
$(\gamma_n)$ be a sequence of hitting times for $A$,
where $A \subset \mathcal S$, then,
for any $N \in \mathbb N$,
and for any $(x_1,S_1) \dots,(x_N,S_N) \in \chi \times \mathcal S$,
\begin{equation}\label{EQ_LEMMA_STRONG_CONDITIONAL_INDEPENDENCE}
\mathbb P \big(  Y^{\gamma_N+1}_{\gamma_1+1} \in S^{N}_1 \mid  X^{\gamma_N+1}_{\gamma_1+1} = x^{N}_1  \big)
= \prod_{k=1}^N \mathbb P \big( Y_{\gamma_k+1} \in S_{k} \mid X_{\gamma_{k}+1}=  x_{{k}}\big).
\end{equation}
\end{lemma1}

\begin{proof} Let $C:= \big(Y^{\gamma_{N-1}+1}_{\gamma_1+1} \in S^{N-1}_1,  X^{\gamma_{N-1}+1}_{\gamma_1+1} = x^{N-1}_1 \big)$
for readability.
\begin{align*}
\left. \right. &\mathbb P \big(   Y^{\gamma_N+1}_{\gamma_1+1} \in S^{N}_1, X^{\gamma_N+1}_{\gamma_1+1} = x^{N}_1  \big)
=\mathbb P \big( Y_{\gamma_N+1} \in S_{N} , \,  X_{\gamma_N+1} = x_{N}  \mid C\big) \mathbb P(C) \\
&\phantom{\qquad\quad}=\mathbb P \big(  Y_{\gamma_N+1} \in S_{N} , \,  X_{\gamma_N+1} = x_{N}  \mid
  X_{\gamma_{N-1}+1} = x_{N-1}\big) \mathbb P(C) \\
&\phantom{\qquad\quad}= \mathbb P \big(  Y_{\gamma_N+1} \in S_{N}  \mid   X_{\gamma_N+1} = x_{N}  \big)
   \mathbb P \big(   X_{\gamma_N+1} = x_{N} \, \mid X_{\gamma_{N-1}+1} = x_{N-1}   \big) \mathbb P(C) \\
 %
&\phantom{\qquad\quad}=\prod_{k=1}^N \mathbb P \big( Y_{\gamma_k+1} \in S_{k} \mid X_{\gamma_{k}+1}=  x_{{k}}\big)
\mathbb P \big( X^{\gamma_{N-1}+1}_{\gamma_1} = x^{N-1}_1 \big),
\end{align*}

\noindent
where the second equality follows
by Remark \ref{REMARK_STRONG_SPLITTING},
the third by Lemma \ref{LEMMA_READ_OUT_RANDOM_TIMES},
 and the last equality follows iterating the procedure and using Remark \ref{REMARK_STRONG_STRONG_MARKOV}.
\end{proof}

\subsection{HMMs and countable mixtures of i.i.d. sequences}
The following fact was used in the proof of Theorem \ref{THEOREM_DHARMADHIKARI_CONTINUOUS_STATE_SPACE}.
If the HMM $(Y_n)$ has an underlying Markov chain with block structured transition probability matrix, with identical rows within blocks, then $(Y_n)$ is a countable mixture of i.i.d. sequences.

Consider a Markov chain $(X_n)$ with values in $\mathcal X$
and transition matrix $P$ as follows
\begin{eqnarray}\label{TRANSITION_MATRIX_MIXT_IID}
P := \left(
\begin{array}{cccccc}
P^1 &  0  & \ldots & \ldots &  0 &  0  \\
0 & P^2 & 0 & \ldots  & 0 &  0  \\
\vdots & \vdots & \vdots & \vdots & \vdots \\
0 & \ldots& \ldots& 0 & P^h &  0 \\
0&  \dots & \ldots & \ldots& \ldots & \ddots
\end{array} \right),
\quad
P^h := \left(
\begin{array}{cccc}
p^h_{c^h_1} & p^h_{c^h_2} & \ldots &  p^h_{c^h_{l_h}} \\
p^h_{c^h_1} & p^h_{c^h_2} & \ldots &  p^h_{c^h_{l_h}} \\
\vdots & \vdots & \vdots & \vdots \\
p^h_{c^h_1} & p^h_{c^h_2} & \ldots &  p^h_{c^h_{l_h}} \\
\end{array} \right),
\end{eqnarray}
with $h \in H$, a countable set. The block $P_h$ has size $l_h$.
Some of the $p^h_c$ can be null.
The Markov chain $(X_n)$ has clearly $H$ recurrence classes,
one for each block, and no transient states.
Let us indicate with $C_h$ the $h$-th recurrence class,
corresponding to the states of the $h$-th block,
set $C_h = \{c^h_1, \dots, c^h_{l_h}\}$, where $l_h$ can be infinite.
Trivially $\mathcal X = \cup_{h \in H} C_h$.
An invariant distribution associated with the $h$-th block is
$\mathbf{p}^h:=(p^h_{c^h_1}, \dots, p^h_{c^h_{l_h}})$, and
for any sequence $\mu_h > 0$ with $\sum_{h\in H}  \mu_h =1$, the vector
\begin{equation}\label{EQ_PI}
\pi=( \mu_1 \mathbf{p}^1, \dots, \mu_h \mathbf{p}^h, \dots )
\end{equation}
is an invariant distribution for $P$.

\begin{lemma1}\label{LEMMA_LUMPING_CONTINUOUS_READ_OUT}
Consider a HMM $(Y_n)$ where the underlying Markov chain $(X_n)$
has transition matrix $P$ as in (\ref{TRANSITION_MATRIX_MIXT_IID}),
invariant measure $\pi$ as in (\ref{EQ_PI}), and assigned read-out distributions $f_x(\bar S)$,
then $(Y_n)$ is a countable mixtures of i.i.d. sequences
where $\widetilde p$ takes values in the set $\{ F_h, \, h \in H  \}$, with
$$
F_h(\bar S )  := p^h_{c^h_1}  f_{c^h_1}(\bar S) + \dots   +p^h_{c^h_{l_h}} f_{c^h_{l_h}}(\bar S),
$$
and $ \mathbb P(\widetilde p =F_h)=\mu_h$.
\end{lemma1}

\begin{proof}
Let us compute the finite distributions of $( Y_n)$. For any $N \in \mathbb N$,
let $S_0, \dots, S_N \in \mathcal S$:

\begin{align*}
\left. \right. &\mathbb P\big( Y_0^N  \in S_0^N \big)
 =\sum_{x_0^N \in \mathcal X} \mathbb P\big(  Y_0^N \in S_0^N, X_0^N=x_0^N \big) \\
&= \sum_{x_0^N \in \mathcal X}  P \big( X_0= x_0 \big)  \prod_{n=0}^N \mathbb  P( Y_n \in  S_n \mid X_n=x_n )
  \prod_{n=1}^N  \mathbb P(X_n=x_n \mid X_{n-1}=x_{n-1}) \\
&=\sum_{x_0^N \in \mathcal X} \pi_{x_0} \,\prod_{n=0}^N  f_{x_n} (S_n)  \,  \prod_{n=1}^N P_{x_{n-1},x_n}  \\
&=  \sum_{x_0^N \in \mathcal X}
\pi_{x_0} f_{x_0} (S_0)   P_{x_0,x_1} f_ {x_1} (S_1)   \dots   P_{x_{N-1},x_{N}}  f_ {x_{N}} (S_N)  \\
&=\sum_{h\in H} \sum_{x_0^N \in \, C_h}
\mu_h p^{h}_{x_0} f_{x_0} (S_0) \,  p^h_{x_1} \,f_{x_1} (S_1)  \dots   p^h_{x_{N}}  f_ {x_{N}} (S_N)\\
&=  \sum_{h \in H} \mu_h \sum_{x_0 \in \, C_h}
 p^{h}_{x_0} f_{x_0} (S_0)  \Big( \sum_{x_1 \in \, C_h} p^h_{x_1} f_{x_1} (S_1)  \dots \big( \sum_{x_{N}\in \, C_h} p^h_{x_{N}}  f_ {x_{N}} (S_N)\big)\Big)\\
&=   \sum_{h \in H} \mu_h F_h (S_0)  \,F_h (S_1) \, \dots F_h (S_N),
\end{align*}
where the second equality follows by the HMM properties,
the fifth equality follows noting that $P_{x_n, x_{n+1}}$ is null for $x_n$ and $x_{n+1}$
in different recurrence classes, and it is equal to $p^h_{x_{n+1}}$ for $x_n$ and $x_{n+1}$
in the same recurrence class $C_{h}$.
The expression above coincides with the representation of countable mixtures of i.i.d. sequences given in (\ref{IID_COUNTABLE_MIXTURE_CONTINUOUS_STATE_SPACE}), thus completing the proof.
\end{proof}


%

\begin{thebibliography}{00}




\bibitem{ALDOUS} Aldous, D.J. (1985)
    \emph{Exchangeability and related topics}
     in {\slshape Ecole d'\'Et\'e de Probabilit\'es de Saint-Flour XIII - 1983}, Lecture Notes in Mathematics 1117, Springer, Berlin

 \bibitem{DEFINETTI_1938} de Finetti, B. (1938)
 \emph{Sur la condition d' equivalence partielle}
 in {\slshape Actualit\'e Scientifiques et Industrielles},
Hermann, Paris, \textbf{739}, 5-18


\bibitem{DHARMADHIKARI} Dharmadhikari, S.W. (1964)
           \emph{Exchangeable processes which are function of stationary Markov chains}
    in {\slshape The Annals of Mathematical Statistics}, \textbf{35}, 429-430

\bibitem {DIACONIS_FREEDMAN_1980} Diaconis, P. and Freedman, D. (1980)
             \emph{de Finetti's theorem for Markov chains}
    in {\slshape The Annals of Probability}, \textbf{8}, 115-130

\bibitem{DIACONIS-FREEDMAN-regularity1}
Diaconis, P. and Freedman, D. (2004)
\emph{The {M}arkov Moment Problem and de {F}inetti's Theorem, {P}art {I} and {P}art {I}{I}}
in { \slshape Mathematische Zeitschrift},
\textbf{247},
183-212


\bibitem{EPIFANI_FORTINI_LADELLI} Epifani, I. and Fortini, S. and Ladelli, L. (2002)
\emph{A characterization for mixtures of semi-Markov processes},
in {\slshape Statistics and Probability Letters}
\textbf{60}, 445-457


\bibitem{FINESSO_PROSDOCIMI_ECC} Finesso, L. and Prosdocimi, C. (2009)
\emph{Partially exchangeable hidden {M}arkov models},
in {\slshape Proceeding of European Control Conference 2009}, 3910--3914


\bibitem{FORTINI_LADELLI_PETRIS_REGAZZINI} Fortini, S., Ladelli, L., Petris, G. and Regazzini, E. (2002)
         \emph{On mixtures of distribution of Markov chains}
         in {\slshape Stochastic Processes and their Applications}, \textbf{100}, 147-165

\bibitem{KALLENBERG_BOOK_2005} Kallenberg, O. (2005)
in {\slshape Probabilistic Symmetries and Invariance Principles}
, Springer


\bibitem{PROSDOCIMI}  Prosdocimi, C. (2010)
\emph{Partial exchangeability and change detection for hidden {M}arkov
models},
in {\slshape {PhD} Dissertation, cycle XXII, University of Padova }

\bibitem{VIDYASAGAR_2011_REALIZATION} Vidyasagar, M. (2011)
\emph{The complete realization problem for hidden Markov models: A survey and some new results}
in {\slshape Mathematics of Control, Signals and Systems},
\textbf{23} (1), 1-65











 \end{thebibliography}
%



\end{document}